# Polynomial and trigonometric splines

Denysiuk V.P. Dr of phys-math. sciences, Professor, Kiev, Ukraine
National Aviation University
kvomden@nau.edu.ua

## Annotation

Classes of simple polynomial and simple trigonometric splines given by Fourier series are considered. It is shown that the class of simple trigonometric splines includes the class of simple polynomial splines. For some parameter values, the polynomial splines coincide with the trigonometric ones; this allows to transfer to such trigonometric splines all the results obtained for polynomial splines. Thus, it was possible to combine two powerful theories - the theory of trigonometric Fourier series and the theory of simple polynomial splines. The above material is illustrated by numerous examples.

**Keywords**:

Generalized trigonometric functions, interpolation, polynomial and trigonometric splines.

## Introduction

Approximation, respectively, the representation of a known or unknown function through a set of some special functions can be considered as a central topic of analysis; such special functions are well defined, easy to calculate, and have certain analytical properties [ 1]. Algebraic and trigonometric polynomials, exponential functions often act as special functions. [2], multinomial [3] and trigonometric [4] splines, etc.

For the whole $r \geq 0$ by $C^r_{[a,b]}$ denote the set $r$ times continuously differentiated on $[a,b]$ functions $r$ times continuously differentiated on $[a,b]$ functions, and the derivative of the order $r+1$ is a function of limited variation. Let us also denote by $C^{-1}_{[a,b]}$ - the set of piecewise continuous functions with breakpoints of the first kind. Also through $C^r_p$ denote the set $2\pi$ - periodic, $r$ times continuously differentiated on $[0, 2\pi)$ functions, and through $C^{-1}_p$ - the plural $2\pi$ periodic, piecewise continuous functions with breakpoints of the first kind. Note that another system of notation is often used.

Known for [5,6], which is the best apparatus for approximating the functions of a class $C^r_p$ there are simple polynomial splines. The theory of such splines is well developed ([3], [7], [8], [9]) etc.. The main disadvantage of polynomial splines, in our opinion, is that they have a lumpy structure; this leads to the fact that in practice mainly third-order splines are used, which are sewn from pieces of algebraic polynomials of the third degree. However, this structure of splines significantly limits their application in many problems of computational mathematics.

(Schoenberg) [10], proposed trigonometric splines that are crosslinked from trigonometric polynomials of a certain degree, and have the same disadvantage as polynomial splines - a lump structure.

In [4], [11], [12] another principle of construction of trigonometric splines was proposed, their representation by uniformly convergent trigonometric series (by Fourier series), whose coefficients have a certain descending order. The undeniable advantage of such splines is that they are given as a single expression over the entire interval of the function.

One of the properties of the trigonometric splines introduced in this way is that at certain values of the parameters they coincide with simple polynomial splines. [11]; therefore, in this case, all the results of the approximation estimates can be transferred to the trigonometric splines, obtained for polynomial splines.

## The purpose of the work.

Construction of classes of interpolation trigonometric splines depending on parameter vectors $\Gamma$ and $H$, and the study of some of their properties.

## The main part.

Let a segment be given $[a,b]$ and let some partition be given on this segment $\Delta_N = \{x_i\}_{i=1}^N$, $a = x_1 < x_1 < ... < x_N = b$, $N$ natural Let also set the sequence of values $y_1, y_2, ..., y_N$.

Consider the problem of interpolation in such a formulation. Find a function $f(x)$, $f(x) \in C_{[a,b]}^r$, which satisfies the conditions $f(x_i) = y_i$, $i = 1, 2, ..., N$.

It is clear that in such a statement of the interpolation problem the main difficulty is to construct $r$ times continuously differentiated by $[a,b]$ functions $f(x)$, depending on $N$ parameters. Such a construction can currently be done in two ways.

With one of them a function $f(x) \in C_{[a,b]}^r$ constructed by the method of piecewise polynomial or piecewise trigonometric interpolation. In the case of piecewise polynomial interpolation, we arrive at a well-developed theory of polynomial splines. (we will indicate only the classic works of Alberg, Stechkin, Zavyalov, Tikhomirov, Korneychuk and many others); in the case of piecewise-trigonometric interpolation we obtain the already mentioned Schoenberg splines [10], whose theory is relatively poorly developed. Note that the attempts of lumpy crosslinking of generalized polynomials composed of functions of other classes are unknown to the author. The main disadvantage of this approach is that as a result of performing differentiation and integration operations again get piecewise-polynomial or piecewise-trigonometric splines.

In another way, the function $f(x) \in C_{[a,b]}^r$ are constructed as the sum of infinite, uniformly convergent Fourier trigonometric series. This method is based on the theory of generalized trigonometric functions [11], which are a generalization of periodic functions (distributions) L. Schwartz [13]. For some values of the parameters that determine these trigonometric series, both methods lead to the same results; this gives grounds to call the functions obtained in this way trigonometric splines. In addition, this fact allows to transfer to such trigonometric splines all the results of the theory of approximations obtained for polynomial splines.

Thus, it was possible to combine two powerful theories - the theory of polynomial splines and the theory of trigonometric Fourier series; this gives reason to expect new results of both theoretical and practical nature.

Note that, generally speaking, the term "trigonometric splines" is often used to refer to Schoenberg splines obtained by the method of piecewise trigonometric interpolation; however, the author, avoiding the introduction of new terms, considers it possible to transfer the term "trigonometric splines" to functions with similar properties obtained through trigonometric series.

Both methods of constructing polynomial and trigonometric splines have their advantages and disadvantages.

In what follows, we will consider only simple polynomial splines, ie splines of the defect 1 [3]. It is clear that a simple degree spline $r$ belongs to $C_{[a,b]}^{r-1}$.

Indexing of trigonometric splines can be entered in two ways. At one of them the order $r$ of a trigonometric spline is determined by the degree of a simple polynomial spline, which is analogous to a trigonometric spline. In another way, the order $r$ the trigonometric spline will determine the number of continuous derivatives of this spline. Although in the vast majority of cases trigonometric splines have no polynomial analogues, still for greater clarity we choose the first method of indexing.

Consider the main properties of polynomial and trigonometric splines, which are related to the methods of their construction, in more detail; the order of coverage of these properties by the author is chosen arbitrarily.

**Periodicity.**

When constructing polynomial splines, there are periodic and non-periodic splines; this difference is taken into account by setting certain boundary conditions. However, using the method of phantom nodes of periodic continuation of functions offered by the author $f(x) \in C_{[a,b]}^k$, non-periodic splines are easy to reduce to periodic [14]. The expediency of such a summary is explained by the fact that the fundamental results of the theory of approximations with respect to polynomial splines are obtained precisely for the periodic case. [5].

Trigonometric splines by construction are periodic and belong to the set $C_{[a,b]}^k$; however, they can also be used to interpolate non-periodic functions.

**Convergence factors.**

When constructing trigonometric splines, it is necessary to specify the type of convergence factors [Den], which have a descending order $O(n^{-(1+r)})$, in which the sum of the trigonometric series belongs to the set $C_{[0,2\pi]}^{r-1}$. Polynomial analogues of trigonometric splines exist, for example, in the case where the convergence factors have the form $\left[\frac{1}{k}\sin\left(\frac{\pi k}{N}\right)\right]^{1+r}$. In the future, we will limit ourselves to considering such a case.

**Degrees of splines.**

When constructing simple polynomial splines, we have to solve systems of linear equations with tape matrices, the complexity of which increases sharply with increasing degrees of splines; this leads to the fact that in practice the most commonly used splines of the third degree.

When constructing trigonometric splines there are no restrictions on their order.

**Parity of degrees.**

When constructing polynomial splines, it is necessary to distinguish between even and odd values of the degrees of the spline, because the algorithms for constructing such splines differ. At construction of trigonometric splines such need does not arise as algorithms of construction of such splines are identical.

**Uniform crosslink grids and interpolation grids.**

When constructing polynomial splines of odd degree, the crosslink grid is also an interpolation grid; when constructing polynomial splines of even degree, the crosslink grid differs from the interpolation grid.

When constructing trigonometric splines, it is necessary to enter the definition of the crosslink grid. The crosslink grid will be called the grid in the nodes of which the polynomial analogue (polynomial spline of the first degree) of the trigonometric spline of the 1st order is stitched..

The crosslink grid and the interpolation grid can be the same or different; this is given by the parameters of the trigonometric spline. It follows that not all trigonometric splines have polynomial analogues.

**Uneven crosslink grids and interpolation grids.**

Construction of polynomial splines on such grids does not cause special difficulties.

The construction of trigonometric splines on non-uniform grids was not considered by the author.

**Construction algorithms.**

As we have already said, when constructing polynomial splines we have to solve systems of linear equations with tape matrices, the complexity of which increases sharply with increasing degree of splines..

Discrete Fourier coefficients have to be calculated when constructing trigonometric splines, and the well-known Fast Fourier Transform (FFT) algorithms can be widely used. When constructing trigonometric splines, there are no restrictions on their order.

**Number of interpolation nodes.**

When constructing polynomial splines, there are no fundamental restrictions on the number of interpolation points.

When constructing trigonometric splines, the restriction on the number of interpolation points is imposed only by FFT algorithms, if they are used. If such algorithms are not used, then there are no fundamental restrictions on the number of interpolation points.

**Fundamental interpolation splines**

The construction of fundamental interpolation polynomial splines is quite complex; such splines have mainly theoretical applications.

Construction of fundamental interpolation trigonometric splines does not cause special difficulties; such splines can be recommended for widespread use in practice.

**Fundamental approximation splines.**

The author does not know about the fundamental approximation polynomial splines.

Fundamental approximation trigonometric splines are constructed naturally, based on the approximation of fundamental trigonometric splines.

**B-splines.**

Polynomial B-splines are constructed sequentially with increasing degree or by convolution, or based on recurrent relations.

Trigonometric B-splines are constructed as the sums of fast-converging trigonometric series; there is no need for their consistent construction with increasing degree.

**Interpolation and apoximation splines of several variables.**

Algorithms for constructing interpolation polynomial splines of several variables are quite complex.

Interpolation trigonometric splines of several variables are constructed as products of fundamental interpolation trigonometric splines on each of variables..

The same applies to approximation splines of several variables.

Dependency introduced trigonometric splines from vectors $\Gamma = \{\gamma_1, \gamma_2, \gamma_3\}$ and $H = \{\eta_1, \eta_2, \eta_3\}$. The splines thus obtained belong to the set $C^k_{[a,b]}$, interpolate a given sequence of values, but have polynomial analogues only in the case of single vectors $\Gamma$ and $H$. interpolate a given sequence of values, but have polynomial analogues only in the case of single vectors.

**Differentiation and integration.**

As we have already said, the main disadvantage of polynomial splines is their piecewise polynomial structure; accordingly, the operations of differentiation and integration of such splines again lead to piecewise polynomial functions.

The operations of differentiation and integration of trigonometric splines are performed much easier, namely by articulated differentiation or integration of uniformly convergent trigonometric series. The results of these operations are easy to obtain by changing the parameter $q$; it is clear that when $0 < q \leq r$ we have the operation of differentiation, and when $q < 0$ - integration operation. Note that for integer values of this parameter we have ordinary operations of differentiation and integration, and for non-integer values we obtain fractional derivatives and fractional integrals in Weyl's understanding. [15]. The use of fractional derivatives and integrals of trigonometric splines is quite interesting, but we have not studied.

As an illustration to the above, consider the construction of a trigonometric spline on $[0, 2\pi)$. We introduce on $[0, 2\pi)$ uniform grids $\Delta_N^{(I)}$, $I$ - grid indicator, ($I = 0,1$), and $\Delta_N^{(0)} = \{x_i^{(0)}\}_{i=1}^N$, $x_i^{(0)} = \frac{2\pi}{N}(i-1)$, and $\Delta_N^{(1)} = \{x_i^{(1)}\}_{i=1}^N$, $x_i^{(1)} = \frac{\pi}{N}(2i-1)$, $N = 2n+1$, $n = 1, 2, ...$. In addition, denote by $a_0, a_k, b_k$ and $a_0, a1_k, b1_k$, ($k = 1, 2, ..., n$) coefficients of interpolation trigonometric polynomials interpolating given values $y_1, y_2, ..., y_N$ in nodes respectively grids $\Delta_N^{(0)}$ and $\Delta_N^{(1)}$, and which are calculated by known formulas (see. [2]).

We present trigonometric interpolation splines on grids in general as follows [12]

$$St(I_1, I_2, \Gamma, H, \nu, r, q, t) = \frac{a_0^{(I_2)}}{2} I(q) + \sum_{k=1}^{\frac{N-1}{2}} \left[ a_k^{1-I_2} a1_k^{I_2} \frac{C_k(I_1, \Gamma, \nu, r, q, t)}{hc_k(I_1, I_2, \Gamma, r)} + b_k^{1-I_2} b1_k^{I_2} \frac{S_k(I_1, H, \nu, r, q, t)}{hs_k(I_1, I_2, H, r, k)} \right],$$

where

$$I(q) = \begin{cases} 1, & \text{if } q = 0; \\ 0, & \text{if } q \neq 0. \end{cases} \text{ - derivative indicator;}$$

$$C_k(I_1, \Gamma, \nu, r, q, t) = \gamma_1 \nu_k(r) k^q \cos(kt + \frac{\pi}{2} q) +$$
$$+ \sum_{m=1}^{\infty} (-1)^{m I_1} \left[ \gamma_3 \nu_{mN+k}(r)(mN+k)^q \cos((mN+k)t + \frac{\pi}{2} q) + \gamma_2 \nu_{mN-k}(r)(mN-k)^q \cos((mN-k)t + \frac{\pi}{2} q) \right];$$

$$S_k(I1, H, \nu, r, q, t) = \eta_1 \nu_k(r) k^q \sin(kt + \frac{\pi}{2} q) +$$
$$+ \sum_{m=1}^{\infty} (-1)^{m I_1} \left[ \eta_3 \nu_{mN+k}(r)(mN+k)^q \sin((mN+k)t + \frac{\pi}{2} q) - \gamma_2 \nu_{mN-k}(r)(mN-k)^q \sin((mN-k)t + \frac{\pi}{2} q) \right];$$

with convergence factors $\nu_k(r)$, having a descending order $O(k^{-(1+r)})$, and interpolation factors

$$hc_k(I_1, I_2, \Gamma, \nu, r) = \gamma_1 \nu_k(r) + \sum_{m=1}^{\infty} (-1)^{m(I_1 - I_2)} \left[ \gamma_3 \nu_{mN+k}(r) + \gamma_2 \nu_{mN-k}(r) \right];$$

$$hs_k(I_1, I_2, H, \nu, r) = \eta_1 \nu_k(r) + \sum_{m=1}^{\infty} (-1)^{m(I_1 - I_2)} \left[ \eta_3 \nu_{mN+k}(r) + \eta_2 \nu_{mN-k}(r) \right]. \tag{4}$$

where the indicator $I_1$ ($I_1 = 0,1$) determines the stitching grid, indicator $I_2$ ($I_2 = 0,1$) determines the interpolation grid, $a_0^{(I_2)}, a_k^{(I_2)}, b_k^{(I_2)}$ - coefficients of interpolation trigonometric polynomial on the grid $\Delta_N^{(I_2)}$, $r$, ($r = 1, 2, ...$) - a parameter that determines the order of the spline, $q$ - a parameter that determines the order of the derivative; $\Gamma = \{\gamma_1, \gamma_2, \gamma_3\}$ and $H = \{\eta_1, \eta_2, \eta_3\}$ - parameter vectors, and parameters $\gamma_k$ and i $\eta_k$, ($k = 1, 2, 3$) take arbitrary valid values, and at least one of the parameters $\gamma_2$, $\gamma_3$, $\eta_2, \eta_3$ not equal to 0, since when $\Gamma = \{1, 0, 0\}$ and $H = \{1, 0, 0\}$ we have a regular trigonometric interpolation polynomial $T_n^{(I_2)}(t)$. Note that to reduce the entries in the symbols of the splines and the functions through which they are built, we omit the dependence on the number of nodes of the interpolation grids $\Delta_N^{(I)}$.

It is easy to see that in the accepted notation $St(I_1, I_2, \Gamma, H, \nu, r, q, t) \in C_p^{r-1-q}$

Trigonometric splines $St(I_1, I_2, \Gamma, H, \nu, r, q, t)$ we will call simple if $\Gamma = \{1,1,1\}$ and $H = \{1,1,1\}$; in the notation of simple trigonometric splines we will omit the dependence on the vectors $\Gamma$ and $H$.

Quite interesting, in our opinion, is the fact that in some cases interpolation factors $hc_k(I_1, I_2, \nu, r)$ and $hs_k(I_1, I_2, \nu, r)$ of simple trigonometric splines to the nearest sign coincide with the coefficients of trigonometric polynomials, which provide the best approximation of Bernoulli functions in the space metricy $L$ [5].

Here is an example that illustrates the behavior of a simple trigonometric spline on different grids of crosslinking and interpolation.

Example. Let's put it down $N = 9$ and set to $[0, 2\pi)$ grids $\Delta_N^{(0)}$ and $\Delta_N^{(1)}$. Let also be a sequence of values 2,1,3,2,4,1,3,1,3. Calculate the coefficients of interpolation trigonometric polynomials $T0_4(t)$ and $T1_4(t)$, interpolating a given sequence in grid nodes $\Delta_N^{(0)}$ and $\Delta_N^{(1)}$ in accordance.

Graphs of interpolation trigonometric polynomials $T0_4(t)$ and $T1_4(t)$ and graphs of trigonometric splines $St(0,0,1,0,t)$, $St(1,0,1,0,t)$, $St(0,1,1,0,t)$ and $St(1,1,1,0,t)$ shown in fig. 1.

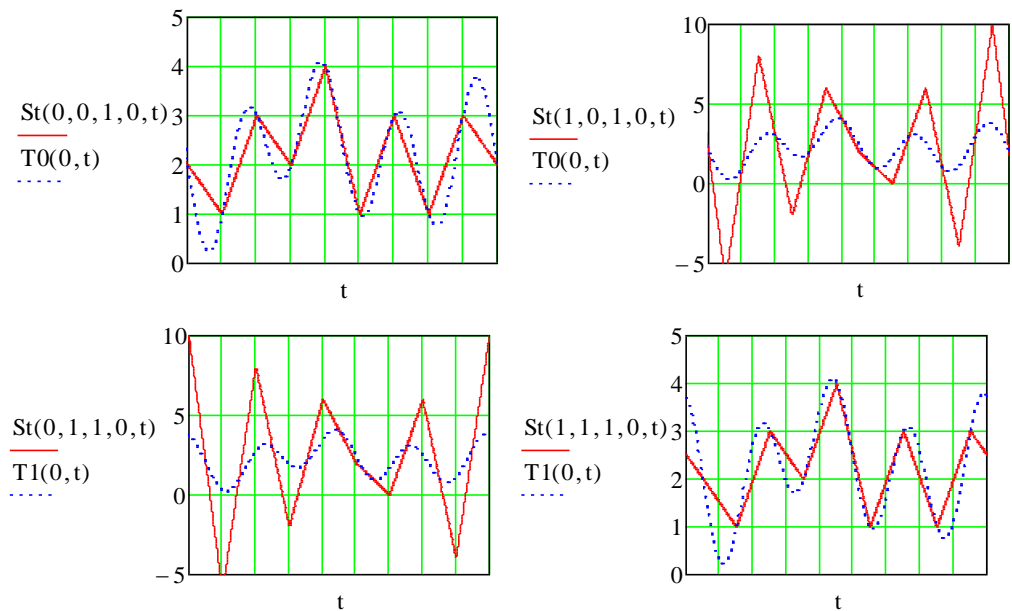

Fig.1. Graphs of trigonometric splines $St(0,0,1,0,t)$, $St(1,0,1,0,t)$, $St(0,1,1,0,t)$ and $St(1,1,1,0,t)$ first order.

It's easy to see that spline $St(0,0,1,0,t)$ is sewn in grid knots $\Delta_N^{(0)}$ and interpolates the specified values in the nodes of the same grid. Spline $St(1,0,1,0,t)$ is sewn in grid knots $\Delta_N^{(1)}$ and interpolates the setpoints in the grid nodes $\Delta_N^{(0)}$. Spline $St(0,1,1,0,t)$ on the contrary, it is sewn in grid knots $\Delta_N^{(0)}$ and interpolates the setpoints in the grid nodes $\Delta_N^{(1)}$. Finally the spline $St(1,1,1,0,t)$ is stitched in grid nodes and interpolates setpoints in nodes of the same grid.

It is clear that of the above splines only spline $St(0,0,1,0,t)$ has a polynomial analogue. Polynomial analogues also have splines $St(0,0,2j+1,0,t)$ ($j = 1, 2, ...$), odd degree, however, the construction of these analogues at $j > 2$ quite complex. It is easy to write polynomial analogues for other types of the given trigonometric splines, but the author does not know about such analogues.

We now show in Fig.2 graphs of splines $St(0,0,2,0,t)$, $St(1,0,2,0,t)$, $St(0,1,2,0,t)$ і $St(1,1,2,0,t)$ парного степеня.

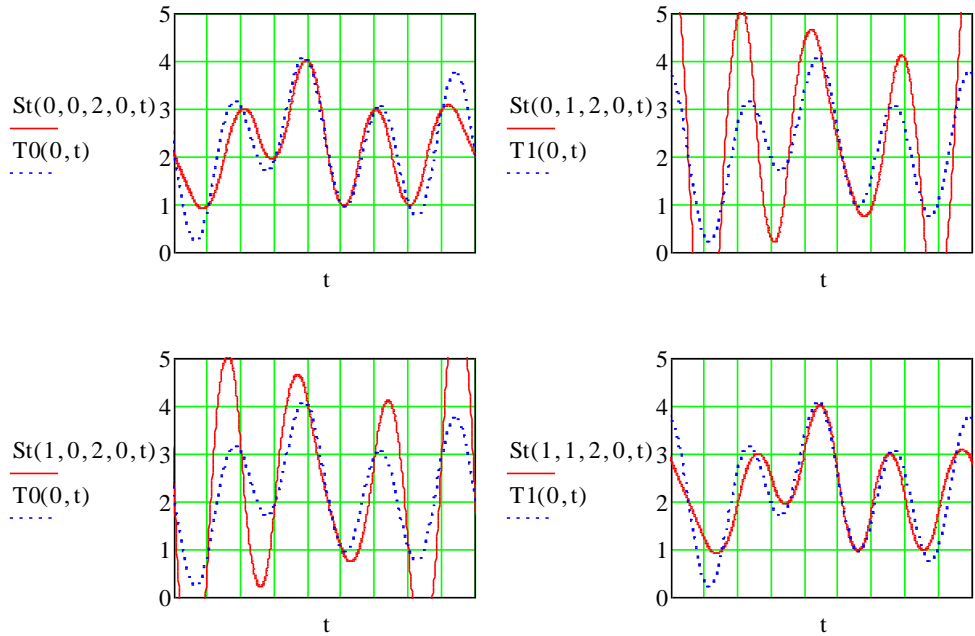

Fig.2. Spline charts $St(0,0,2,0,t)$, $St(1,0,2,0,t)$, $St(0,1,2,0,t)$ and $St(1,1,2,0,t)$ even degree.

It's easy to see that spline $St(0,0,2,0,t)$ is sewn in grid knots $\Delta_N^{(0)}$ and interpolates the specified values in the nodes of the same grid. Spline $St(0,1,2,0,t)$ is sewn in grid knots $\Delta_N^{(0)}$ and interpolates the specified values in the nodes of the Spline grid $St(1,0,2,0,t)$ on the contrary, it is sewn in grid knots $\Delta_N^{(1)}$ and interpolates the setpoints in the grid nodes $\Delta_N^{(0)}$ $\Delta_N^{(1)}$. Finally the spline $St(1,1,2,0,t)$ is sewn in grid knots $\Delta_N^{(1)}$ and interpolates the specified values in the nodes of the same grid.

It is clear that of the above splines only spline $St(0,1,1,0,t)$ has a polynomial analogue. Polynomial analogues also have splines $St(0,0,2j,0,t)$ ($j = 1,2,...$) even degree, but about the construction of such analogues at $j \geq 2$ the author is unknown.

Finally, Figure 3 shows the graphs of the trigonometric spline, for exampleад, $St(0,0,r,0,t)$ in ascending order of the parameter $r$.

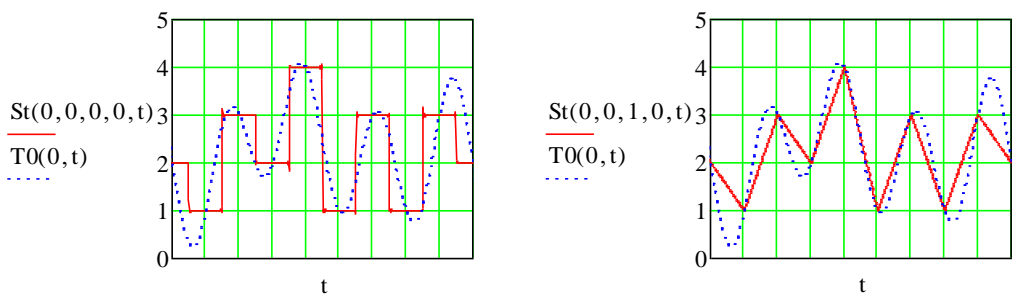

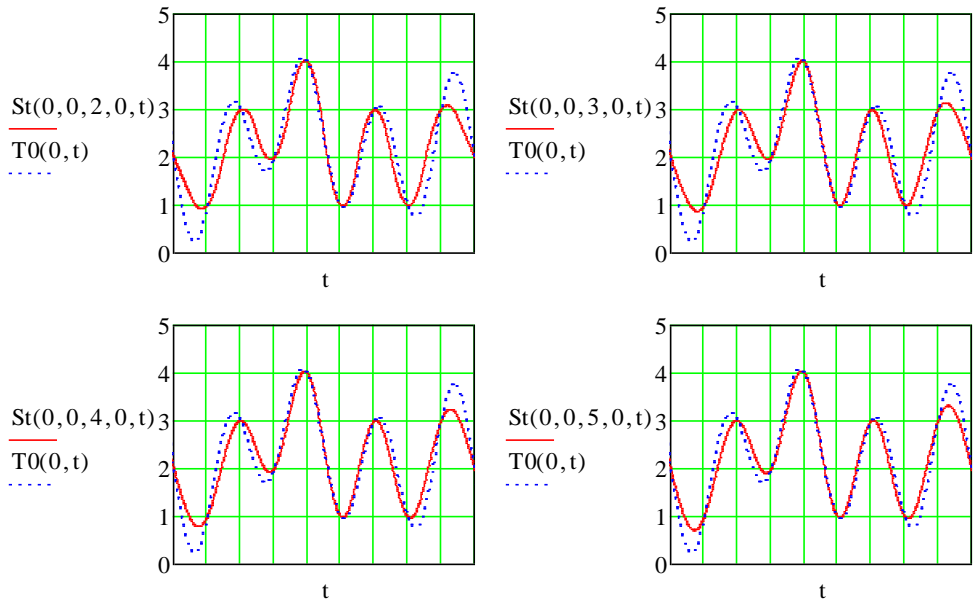

Figure 3. Graphs of the trigonometric spline $St(0,0,r,0,t)$
in ascending order of parameter $r$.

These graphs illustrate the fact that $St(0,0,r,0,t) \in C_p^{r-1}$ for all integer values of the parameter $r$. Also note that the splines are of odd order $St(0,0,1,0,t)$, $St(0,0,3,0,t)$, $St(0,0,5,0,t)$ have polynomial analogues - simple polynomial splines of the 1st, 3rd and 5th degrees, respectively. The polynomial analogue also has a spline of the 0th order - piecewise - steel, and in this case the uniform convergence of trigonometric Fourier series is lost. Splines of the same order $St(0,0,2,0,t)$ i $St(0,0,4,0,t)$ have no polynomial analogues.

Finally, to illustrate, we present graphs of first-order trigonometric splines with vectors $\Gamma = \{1.5, .5, -1\}$ and $H = \{1, .5, .5\}$.

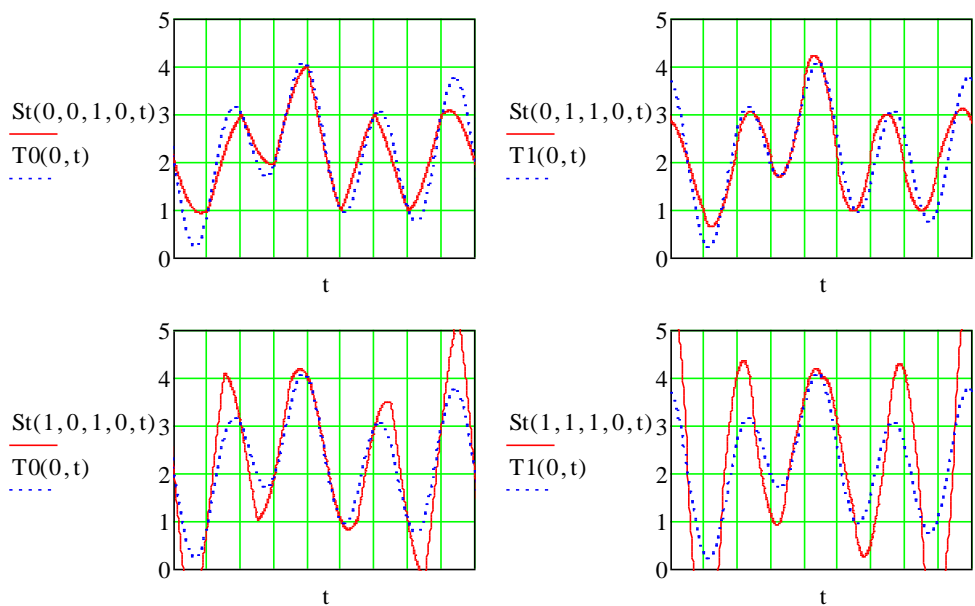

Fig.4. Graphs of first-order trigonometric splines with vectors
$\Gamma = \{1.5, .5, -1\}$ and $H = \{1, .5, .5\}$.

It is easy to see that the given trigonometric splines retain differential and interpolation properties. However, their behavior differs from the behavior of simple first-order trigonometric splines with single vectors shown in Fig.1.

Finally, we formulate some interesting, in our opinion, questions.

1. It is clear that trigonometric splines $St(1,1,r,0,t)$ and $St(1,0,r,0,t)$ are shifted half-step grid splines $St(0,0,r,0,t)$ and $St(1,0,r,0,t)$. Because splines $St(0,0,2j+1,0,t)$ and $St(0,1,2j,0,t)$ ($j=1,2,...$) have polynomial analogues, they are covered by all the results of the theory of approximations obtained for simple polynomial splines. Is it possible to extend the same results to splines $St(1,1,2j+1,0,t)$ and $St(1,0,2j,0,t)$?
2. In the general case, trigonometric splines depend on the parameter vectors $\Gamma$ and $H$. It is clear that any functionalities from trigonometric splines or their derivatives will also depend on these vectors; the study of such dependences can lead to new results in the theory of approximations.

## Conclusions.

1. It is shown that the class of trigonometric splines is represented by Fourier series.
2. Thus it was possible to combine two powerful theories - the theory of polynomial splines and the theory of trigonometric Fourier series; this gives reason to expect new results of both theoretical and practical nature.
3. The set of simple polynomial splines is a subset of the set of trigonometric splines.
4. At some values of parameters trigonometric splines coincide with simple polynomial splines; on such trigonometric splines it is possible to transfer all results of the theory of approximations received for simple polynomial splines.
5. FFT algorithms are widely used in the construction of trigonometric splines.
6. Dependences of trigonometric splines on vectors are introduced $\Gamma$ and $H$. The study of such dependencies looks very promising.
7. Undoubtedly, further studies of trigonometric splines are of undeniable interest.

## List of references